\documentclass[12pt]{amsart}
\usepackage{amscd}     
\usepackage{amssymb}
\usepackage{amsmath, amsthm, graphics}
\usepackage{xypic}     
\LaTeXdiagrams        
\usepackage[all]{xy}
\xyoption{2cell} \UseAllTwocells \xyoption{frame} \CompileMatrices
\allowdisplaybreaks[3]
\usepackage{amsfonts}
\usepackage[normalem]{ulem}
\usepackage{hyperref}
\usepackage{tikz}
\usepackage{mathtools}
\usepackage{verbatim} 
\usepackage{bbold}

\usepackage{latexsym}
\usepackage{epsfig}
\usepackage{amsfonts}
\usepackage{enumerate}
\usepackage{times}

\newtheorem{theorem}{Theorem}[section]
\newtheorem{corollary}[theorem]{Corollary}
\newtheorem{lemma}[theorem]{Lemma}
\newtheorem{example}[theorem]{Example}

\newtheorem{proposition}[theorem]{Proposition}

\newtheorem{definition}[theorem]{Definition}

\theoremstyle{remark}

\theoremstyle{remark}

\numberwithin{equation}{section}

\newcommand{\Mfrak}{\mathfrak{M}}
\newcommand{\tCH}{\mathrm{CH}}

\def\<{\left\langle}
\def\>{\right\rangle}

\def\b1{{\mathbf 1}}

\begin{document}

\title{On the tautological rings of stacks of twisted curves}

\author[Tseng]{Hsian-Hua Tseng}
\address{Department of Mathematics\\ Ohio State University\\ 100 Math Tower, 231 West 18th Ave. \\ Columbus,  OH 43210\\ USA}
\email{hhtseng@math.ohio-state.edu}

\date{\today}



\maketitle

\begin{abstract}
We introduce the tautological rings of moduli stacks of twisted curves and establish some basic properties.
\end{abstract}

\setcounter{tocdepth}{1}

\setcounter{section}{-1}

\section{Introduction}
We work over an algebraically closed field of characteristic $0$.

\subsection{Backgroud}
Balanced twisted curves\footnote{Twisted curves are also called stacky curves or orbifold curves.} are nodal curves with specific types of stack structures possibly at nodes and isolated smooth points. \'Etale locally at a smooth stacky point, it looks like $[\mathbb{A}^1/\mu_r]$ for some $r$, where $\mu_r$ acts on $\mathbb{A}^1$ by multiplications. \'Etale locally at a stacky node, it looks like $[(xy=0)/\mu_s]$ for some $s$, where $\mu_s$ acts\footnote{This particular action is called {\em balanced}.} via $\zeta\cdot (x,y)=(\zeta\cdot x, \zeta^{-1}\cdot y)$.

Twisted curves arise naturally in \cite{av}, \cite{acv} as a solution to the problem of compactifying spaces of maps from nodal curves to a proper Deligne-Mumford stack. More precisely, to compactify a space of maps from smooth curves of a proper Deligne-Mumford stack, one must add maps from nodal twisted curves, not just maps from nodal curves.

Later on, compact moduli spaces of certain maps from balanced\footnote{The balanced condition is used in certain axioms in Gromov-Witten theory.} twisted curves to a smooth proper Deligne-Mumford stack were used to construct Gromov-Witten theory of Deligne-Mumford stack, see \cite{agv}. 

Let 
\begin{equation*}
    \mathfrak{M}_{g,n}^{tw}
\end{equation*}
be the stack of prestable balanced twisted curves of genus $g$ with $n$ marked gerbes. The existence of $\Mfrak_{g,n}^{tw}$ and some basic properties have been known to the experts. It is known that $\Mfrak_{g,n}^{tw}$ is a smooth Artin stack locally of finite type. A construction of $\Mfrak_{g,n}^{tw}$ using logarithmic geometry, and proofs of some of its basic properties, can be found in \cite{ol}. 

A variant of $\Mfrak_{g,n}^{tw}$ was introduced in \cite[Section 2]{c} for the purpose of studying relations between higher genus Gromov-Witten theory and genus $0$ Gromov-Witten theory of symmetric products. This variant was also used in \cite{ajt1} to study genus $0$ Gromov-Witten theory of root gerbes. In this note, this variant is denoted by $\Mfrak_{g,I,m,a}^{tw, triv}$, see (\ref{eqn:moduli_tw_curve_w_label}) below.

\subsection{This paper}
In the study of Gromov-Witten theory of Deligne-Mumford stacks, the stack $\Mfrak_{g,n}^{tw}$ is mostly used as a placeholder. On the other hand, $\Mfrak_{g,n}^{tw}$ is a natural generalization of both the stack $\overline{\mathcal{M}}_{g,n}$ of stable curves and the stack $\mathfrak{M}_{g,n}$ of pretable curves. $\Mfrak_{g,n}^{tw}$ has rich geometry and deserves to be studied further. Motivated by this, this note considers one geometric aspect of $\Mfrak_{g,n}^{tw}$, namely {\em the tautological ring} of $\Mfrak_{g,n}^{tw}$. The approach taken here closely follows the work \cite{bs}. 

Section \ref{sec:def} presents a definition of the tautological rings of $\Mfrak_{g,n}^{tw}$ as a subring of the Chow rings of $\Mfrak_{g,n}^{tw}$. In this construction, the aforementioned variant $\Mfrak_{g,I,m,a}^{tw, triv}$ introduced in \cite[Section 2]{c} plays an important role.  A key reason for working with $\Mfrak_{g,I,m,a}^{tw, triv}$ is that their universal curves have modular descriptions, given by the morphisms of forgetting a non-stacky marked point. This allows the tautological rings of $\Mfrak_{g,I,m,a}^{tw, triv}$ be defined using natural gluing and forgetful maps, see Definition \ref{def:taut_ring} for more details. The tautological rings of $\Mfrak_{g,n}^{tw}$ are then defined using tautological rings of $\Mfrak_{g,I,m,a}^{tw, triv}$ via restrictions, see Definition \ref{def:taut_ring2} for a precise description.

A natural class of elements in the tautological rings, namely {\em decorated strata classes}, are also introduced, see Definition \ref{def:deco_stra_class}.

Section \ref{sec:basis} presents an additive spanning set of the tautological rings. Namely, it is shown that the decorated strata classes span the tautological rings as $\mathbb{Q}$-vector spaces. See Theorem \ref{thm:Q-span} for the precise result.

\subsection{Outlook}
Naturally, the next goal in this direction is to understand the genus $0$ situation, thus extending the work \cite{bs2} to twisted curves. More precisely, it is expected that the Chow ring and tautological ring of $\Mfrak_{0,n}^{tw}$ coincide, and relations in the tautological ring of $\Mfrak_{0,n}^{tw}$ are all obtained from WDVV relations. This will be discussed elsewhere.

A long-term, more ambitious goal, is to understand {\em relations} in the tautological rings of $\Mfrak_{g,n}^{tw}$. 

\subsection{Acknowledgment}
The author thanks J. Schmitt for valuable discussions. The author also thanks the referees for helpful suggestions. This work is supported in part by Simons Foundation Collaboration Grant.

\section{Definitions of tautological rings}\label{sec:def}
There are $n$ universal cyclic gerbes over $\Mfrak_{g,n}^{tw}$ corresponding to marked points. Their fiber product over $\Mfrak_{g,n}^{tw}$, denoted by $\Mfrak_{g,n}^{tw, triv}$, parametrizes twisted curves with {\em trivialized} marked gerbes. The map 
\begin{equation}\label{eqn:univ_gerbe}
    \Mfrak_{g,n}^{tw, triv}\to \Mfrak_{g,n}^{tw}
\end{equation}
is finite and \'etale. 

\subsection{Chow rings}
By \cite[Corollary 1.11]{ol}, $\Mfrak_{g,n}^{tw}$ and $\Mfrak_{g,n}^{tw,triv}$ are locally of finite type. The construction of \cite[Appendix A]{bs} applies to define the Chow groups\footnote{Chow groups are taken with $\mathbb{Q}$-coefficients.} 
\begin{equation*}
    \tCH_*(\Mfrak_{g,n}^{tw}), \quad \tCH_*(\Mfrak_{g,n}^{tw, triv}),
\end{equation*}
which admit usual properties of Chow groups of finite type stacks. By the analysis of stabilizers of twisted curves in \cite{acv}, $\Mfrak_{g,n}^{tw}$ and $\Mfrak_{g,n}^{tw, triv}$ have affine stabilizer groups at geometric points (except when $(g,n)=(1,0)$). Hence intersection products can be defined, making these Chow groups into $\mathbb{Q}$-algebras.

The stack $\Mfrak_{g,n}^{tw}$ is not connected. To specify component(s) of $\Mfrak_{g,n}^{tw}$ to work with, it is convenient to rephrase the setup by introducing labels for marked points. Let $I$ be a set of $n$ elements. Let 
\begin{equation*}
    \Mfrak_{g,I}^{tw}
\end{equation*}
be the stack of prestable balanced twisted curves of genus $g$ and marked gerbes labelled by $I$. Clearly there is an isomorphism $\Mfrak_{g,I}^{tw}\simeq \Mfrak_{g,n}^{tw}$ given by choosing an identification $I\to \{1, 2, ..., n\}$. For a function $m: I\to \mathbb{Z}_{>0}$, let 
\begin{equation*}
    \Mfrak_{g,I,m}^{tw}\subset \Mfrak_{g,I}^{tw}
\end{equation*}
be the locus parametrizing twisted curves whose stack structures at marked points are given by the function $m$. The stacks 
\begin{equation*}
    \Mfrak_{g,I}^{tw, triv}, \quad \Mfrak_{g,I,m}^{tw, triv}
\end{equation*}
are defined by inverse images with respect to (\ref{eqn:univ_gerbe}).

The stacks $\Mfrak_{g,I,m}^{tw}$ and $\Mfrak_{g,I,m}^{tw, triv}$ are smooth Artin stacks locally of finite type. It can be seen that $\Mfrak_{g,I,m}^{tw}$ and $\Mfrak_{g,I,m}^{tw, triv}$ are equidimensional. Furthermore, $\Mfrak_{g,I,m}^{tw}$ and $\Mfrak_{g,I,m}^{tw, triv}$ have good filtrations by finite type substacks, in the sense of \cite[Definition A.2]{bs}, given by loci of twisted curves with bounds on the numbers of nodes.

We are interested in the Chow rings
\begin{equation*}
    \tCH^*(\Mfrak_{g,I,m}^{tw}), \quad \tCH^*(\Mfrak_{g,I,m}^{tw, triv}).
\end{equation*}
Since (\ref{eqn:univ_gerbe}) is a fiber product of cyclic gerbes, with $\mathbb{Q}$-coefficients there is an isomorphism 
\begin{equation*}
   \tCH^*(\Mfrak_{g,I,m}^{tw})\simeq \tCH^*(\Mfrak_{g,I,m}^{tw, triv}) 
\end{equation*}
induced by pulling back via (\ref{eqn:univ_gerbe}). In what follows, we work with $\Mfrak_{g,I,m}^{tw, triv}$.

\subsection{Tautological rings}
In the definition of tautological rings for moduli spaces of {\em stable} curves, the modular interpretation of their universal curves plays an important role. For moduli stacks of {\em prestable} curves, such a moduli interpretation is not valid. This issue is overcome in \cite{bs}. We follow their recipe to treat the case of twisted curves, as follows.

Let $\mathcal{A}$ be a commutative semigroup with unit $\mathbb{0}\in\mathcal{A}$ such that
\begin{enumerate}
    \item 
    $\mathcal{A}$ has indecomposable zero: for $x, y\in \mathcal{A}$, if $x+y=\mathbb{0}$, then we must have $x=y=\mathbb{0}$;
    \item
    $\mathcal{A}$ has finite decomposition: if $a\in \mathcal{A}$, then the set $\{(a_1, a_2)| a_1+a_2=a\}\subset \mathcal{A}\times \mathcal{A}$ is finite.
\end{enumerate}

\begin{example}\label{ex:A0}
The example of $\mathcal{A}$ that matters to us the most is the following: $\mathcal{A}_\mathbb{0}=\{\mathbb{0}, \mathbb{1}\}$ with $\mathbb{0}+\mathbb{0}=\mathbb{0}, \mathbb{0}+\mathbb{1}=\mathbb{1}+\mathbb{0}=\mathbb{1}+\mathbb{1}=\mathbb{1}$.
\end{example}

Fix an $\mathcal{A}$. Let $a\in \mathcal{A}$. We consider the stack 
\begin{equation}\label{eqn:moduli_tw_curve_w_label}
    \Mfrak_{g,I,m,a}^{tw, triv}
\end{equation}
of twisted curves with labellings of irreducible components by elements of $\mathcal{A}$. This stack is defined in \cite[Section 2]{c}, where it is denoted by $\Mfrak_{g,I,m,a}$. Roughly speaking, an object of $\Mfrak_{g,I,m,a}^{tw, triv}$ is a twisted curve $(\mathcal{C}, \{p_i\}_{i\in I})\in \Mfrak_{g,I,m}^{tw, triv}$ together with a map $\mathcal{C}_v\mapsto a_{\mathcal{C}_v}$ from the set of irreducible components $\{\mathcal{C}_v\}$ of the normalization of $\mathcal{C}$ to $\mathcal{A}$ such that $a_{\mathcal{C}_v}$ sum up to $a$. The following stability condition is required: for each irreducible component $\mathcal{C}_v$, either $a_{\mathcal{C}_v}\neq \mathbb{0}$, or $\mathcal{C}_v$ together with its special points is a stable twisted curve.

There is a natural morphism $$    \Mfrak_{g,I,m,a}^{tw, triv}\to     \Mfrak_{g,I,m}^{tw, triv}$$ that forgets the $\mathcal{A}$-valuation.

It follows from \cite[Proposition 2.0.2]{c} and properties of $\Mfrak_{g,I,m}^{tw, triv}$ that  $\Mfrak_{g,I,m,a}^{tw, triv}$ is a smooth Artin stack locally of finite type, with good filtrations by finite type substacks.

Consider the stack $\Mfrak_{g,I\cup\{\bullet\}, m',a}^{tw, triv}$ where $m':I\cup\{\bullet\}\to \mathbb{Z}_{>0}$ is given by $m'|_I=m$ and $m'(\bullet)=1$ (i.e. the marked point $\bullet$ has trivial stack structure). Then there is a morphism 
\begin{equation}\label{eqn:forgetful}
  \pi_{\bullet}:  \Mfrak_{g,I\cup\{\bullet\},m', a}^{tw, triv}\to \Mfrak_{g, I, m,a}^{tw, triv},
\end{equation}
given by forgetting the marked point $\bullet$ and contracting resulting unstable components with $\mathcal{A}$-valuations $\mathbb{0}$. It follows from \cite[Proposition 2.1.1]{c} that (\ref{eqn:forgetful}) is the universal curve.

We now recall the natural gluing maps with target  $\Mfrak_{g,I,m,a}^{tw, triv}$. Let 
\begin{equation}\label{eqn:prestable_graph}
    \Gamma=(V, H, E, L, g:V\to\mathbb{Z}_{\geq 0}, v: H\to V, \iota: H\to H)
\end{equation}
be a prestable graph of genus $g$ with $n$ markings. The required properties of $\Gamma$ are 
\begin{enumerate}
    \item 
    $V$ is the set of vertices, $g:V\to \mathbb{Z}_{\geq 0}$ is the genus assignment;
    \item
    $H$ is the set of half-edges, $v:H\to V$ is the vertex assignment, and $\iota: H\to H$ is an involution;
    \item
    the set of $2$-cycles of $\iota$ in $H$ is by definition the set of edges $E$ (self-edges at vertices are permitted);
    \item
    the set of fixed points of $\iota$ is by definition the set of legs $L$ (which correspond to the $n$ markings);
    \item
    the pair $(V,E)$ defines a connected graph $\Gamma$ satisfying the genus condition 
    \begin{equation*}
        \sum_{v\in V}g(v)+h^1(\Gamma)=g.
    \end{equation*}
\end{enumerate}

For a prestable graph $\Gamma$, we fix an identification $L\simeq I$ between the set of legs and marking labels. In addition, we need to choose a map 
\begin{equation}\label{eqn:orb_str}
    m_\Gamma: H\to \mathbb{Z}_{\geq 0}
\end{equation}
 so that $m_\Gamma|_{L}=m$ and $$m_\Gamma(\iota(h))=m_\Gamma(h)$$ for $h\in H$, which ensures that stack structures at nodes are balanced. We also need to choose an $\mathcal{A}$-valuation $a: V(\Gamma)\to \mathcal{A}$ satisfying $\sum_{v\in V(\Gamma)}a(v)=a$.

For a vertex $v\in V(\Gamma)$, let $H(v)\subset H$ be the set of half-edges at $v$ and $L(v)\subset H(v)$ the set of legs at $v$. $I(v)\subset I$ is defined by the identification $L\simeq I$.

Associated to $(\Gamma, m_\Gamma, a)$ is a gluing morphism
\begin{equation}\label{eqn:gluing_map}
    \xi_{(\Gamma, m_\Gamma, a)}: \Mfrak_{(\Gamma,m_\Gamma,a)}^{tw, triv}=\prod_{v\in V(\Gamma)} \Mfrak_{g(v), H(v),m_\Gamma|_{H(v)},a(v)}^{tw, triv}\to \Mfrak_{g,I,m,a}^{tw, triv},
\end{equation}
which is defined as follows. $\xi_{(\Gamma, m_\Gamma, a)}$ sends a collection $(\mathcal{C}_v, \{p_h\}_{h\in H(v)})_{v\in V}$ to the curve $(\mathcal{C}, \{p_i\}_{i\in I})$ obtained by identifying the markings $p_h, p_{\iota(h)}$ for each pair $(h,\iota(h))$ of half-edges forming an edge of $\Gamma$. Here the gerbe structures at $p_h$ and $p_{\iota(h)}$ are identified via the homomorphism $\mu_{m_\Gamma(h)}\to \mu_{m_\Gamma(\iota(h))}$, $\zeta\mapsto \zeta^{-1}$. Since gerbes at nodes are not trivialized, over its image, $\xi_{(\Gamma, m_\Gamma, a)}$ is the product of universal gerbes at nodes (which correspond to edges of $\Gamma$) and has degree $\prod_{(h, \iota(h))\in E}\frac{1}{m_\Gamma(h)}$, see \cite[Section 5.2]{agv0}.

The following definition is a direct extension of \cite[Definition 1.3]{bs}, which is modelled on the treatment of the tautological rings of moduli stacks of stable curves in \cite{gp03}.

\begin{definition}\label{def:taut_ring}
The tautological rings $$\{\mathrm{R}^*(\Mfrak_{g,I,m,a}^{tw, triv})\}_{g,I,m,a}$$ is the smallest system of unital $\mathbb{Q}$-subalgebras of the Chow rings $\{\tCH^*(\Mfrak_{g,I,m,a}^{tw, triv})\}_{g,I,m,a}$ which is closed under taking pushforwards by the natural forgetful maps (\ref{eqn:forgetful}) and gluing maps (\ref{eqn:gluing_map}).
\end{definition}

Consider the semigroup $\mathcal{A}_\mathbb{0}$ in Example \ref{ex:A0}. Let $\mathbb{1}\in \mathcal{A}_\mathbb{0}$. Consider the subset
\begin{equation*}
    \mathfrak{Z}_{g,I,m}^{tw, triv}\subset \Mfrak_{g,I,m,\mathbb{1}}^{tw, triv}
\end{equation*}
consisting of $\mathcal{A}$-valued curves $(\mathcal{C},\{p_i\}_{i\in I},(a_{\mathcal{C}_v})_v)$ such that one of the values $a_{\mathcal{C}_v}$ is equal to $\mathbb{0}$. The following Lemma is a direct adaptation of \cite[Proposition 2.6]{bs}, whose proof can be easily adapted here as well.

\begin{lemma}\label{lem:open_substack}
$\mathfrak{Z}_{g,I,m}^{tw,triv}$ is closed. The composition 
$$\mathfrak{U}_{g,I,m}^{tw,triv}=\Mfrak_{g,I,m,\mathbb{1}}^{tw, triv}\setminus \mathfrak{Z}_{g,I,m}^{tw,triv}\hookrightarrow \Mfrak_{g,I,m,\mathbb{1}}^{tw, triv}\to \Mfrak_{g,I,m}^{tw, triv}$$
of the inclusion and the natural map defines an isomorphism $$\mathfrak{U}_{g,I,m}^{tw, triv}\simeq \Mfrak_{g,I,m}^{tw, triv}.$$
\end{lemma}

The following extends \cite[Corollary 2.7]{bs}.

\begin{corollary}
The universal curve over $\Mfrak_{g,I,m}^{tw, triv}$ is given by 
\begin{equation*}
    \pi_\bullet: \Mfrak_{g,I\cup\{\bullet\},m, \mathbb{1}}^{tw, triv}\setminus \pi_{\bullet}^{-1}(\mathfrak{Z}_{g,I,m}^{tw,triv})\to \Mfrak_{g, I, m,\mathbb{1}}^{tw, triv}\setminus \mathfrak{Z}_{g,I,m}^{tw,triv}
\end{equation*}
\end{corollary}

\begin{definition}\label{def:taut_ring2}
The tautological ring 
$$\mathrm{R}^*(\Mfrak_{g,I,m}^{tw, triv})\subset \tCH^*(\Mfrak_{g,I,m}^{tw, triv})$$
is defined to be the image of the restriction of $\mathrm{R}^*(\Mfrak_{g,I,m, \mathbb{1}}^{tw, triv})$ to the open substack  from Lemma \ref{lem:open_substack}.
\end{definition}

\subsection{Stratum classes}
We define some basic classes in the tautological ring $\mathrm{R}^*(\Mfrak_{g,I,m,a}^{tw, triv})$. Recall the universal curve $\pi_\bullet$ in (\ref{eqn:forgetful}). Let 
\begin{equation}
    \sigma_i: \Mfrak_{g,I, m, a}^{tw,triv}\to \Mfrak_{g, I\cup\{\bullet\}, m', a}^{tw, triv}, \quad i\in I
\end{equation}
be the section corresponding to the marking labelled by $i\in I$. Note that $\sigma_i$ exists because marked gerbes are trivialized. Let $\omega_{\pi_\bullet}$ be the relative dualizing sheaf of $\pi_\bullet$. 

\begin{definition}\label{def:psi_kappa}
\begin{equation*}
\begin{split}
    &\psi_i:=\sigma_i^*c_1(\omega_{\pi_\bullet})\in \tCH^1(\Mfrak_{g,I, m, a}^{tw,triv}),\quad i\in I\\
    &\kappa_m:=(\pi_\bullet)_*(\psi_\bullet^{m+1})\in \tCH^m(\Mfrak_{g,I, m, a}^{tw,triv}).
\end{split}
\end{equation*}
\end{definition}

\begin{definition}[Decorated stratum classes]\label{def:deco_stra_class}
Let $\Gamma$ be a prestable graph of genus $g$ with $n$ markings. Let $m_\Gamma$ be as in (\ref{eqn:orb_str}) and $a: V(\Gamma)\to \mathcal{A}$ an $\mathcal{A}$-valuation with total value $a\in \mathcal{A}$. A decoration on $(\Gamma, m_\Gamma, a)$ is a class on $\Mfrak_{(\Gamma,m_\Gamma,a)}^{tw, triv}=\prod_{v\in V(\Gamma)} \Mfrak_{g(v), H(v),m_\Gamma|_{H(v)},a(v)}^{tw, triv}$ of the form
\begin{equation}\label{eqn:deco}
    \alpha=\prod_{v\in V}\left(\prod_{i\in H(v)} \psi_{v,i}^{a_i}\prod_{j\geq 1}\kappa_{v, j}^{b_{v,j}} \right)\in \tCH^*(\Mfrak_{(\Gamma,m_\Gamma,a)}^{tw, triv}).
\end{equation}
Here $\psi_{v,i}, \kappa_{v,j}\in \tCH^*(\Mfrak_{g(v), H(v),m_\Gamma|_{H(v)},a(v)}^{tw, triv})$ are $\psi$ and $\kappa$ classes of the moduli space corresponding to the vertex $v\in V(\Gamma)$. The decorated stratum class $[(\Gamma, m_\Gamma, a), \alpha]$ is defined to be the pushforward 
\begin{equation}\label{eqn:deco_class}
    [(\Gamma, m_\Gamma, a), \alpha]:=(\xi_{(\Gamma, m_\Gamma, a)})_*\alpha\in \tCH^*(\Mfrak_{g,I,m,a}^{tw, triv}).
\end{equation}
\end{definition}

\section{An additive spanning set of tautological rings}\label{sec:basis}
The main goal of this section is to prove the following extension of \cite[Theorem 1.4]{bs} and \cite[Proposition 11]{gp03}:
\begin{theorem}\label{thm:Q-span}
The tautological ring $\mathrm{R}^*(\Mfrak_{g,I,m,a}^{tw, triv})$ is spanned by decorated stratum classes (\ref{eqn:deco_class}) as a $\mathbb{Q}$-vector space.
\end{theorem}
The treatment here follows those of \cite{bs} and \cite{gp03}.

Because decorated stratum classes are tautological, their span is a vector subspace of the tautological ring. Since the tautological rings are defined to be the smallest system of unital $\mathbb{Q}$-subalgebras which is closed under pushforwards by (\ref{eqn:forgetful}) and (\ref{eqn:gluing_map}), we must prove two things:
\begin{enumerate}
    \item 
    The span of decorated stratum classes is a $\mathbb{Q}$-subalgebra, i.e. it is closed under intersection product.
    \item
    The span of decorated stratum classes is closed under pushforwards by (\ref{eqn:forgetful}) and (\ref{eqn:gluing_map}), hence the span must be the smallest such system.
\end{enumerate}

The proof involves descriptions of pullbacks and pushforwards of decorated stratum classes under natural maps. We spell out the details in our case, following the treatment in \cite[Section 3.2]{bs}.

\subsection{Closure under intersection product}
We recall the combinatorial set-up. Let $A$ be an $\mathcal{A}$-valued prestable graph. 
\begin{definition}(\cite[Definition 3.5]{bs})
An $A$-structure on an $\mathcal{A}$-valued prestable graph $\Gamma$, denoted by $$\Gamma\to A,$$ consists of choices of subgraphs $\{\Gamma_v\}_{v\in V(A)}$ and maps $V(\Gamma)\to V(A)$ and $H(A)\to H(\Gamma)$ satisfying
\begin{enumerate}
    \item
    The total $\mathcal{A}$-value of $\Gamma_v$ is $a(v)$, where $a : V(A)\to \mathcal{A}$ is the $\mathcal{A}$-valuation of $A$.
    \item 
    The map $V(\Gamma)\to V(A)$ is surjective, and the inverse image of $v\in V(A)$ are precisely the vertices of $\Gamma_v$.
    \item
    The $\mathcal{A}$-valuation of $\Gamma$ is given by the composition of $V(\Gamma)\to V(A)$ and $a: V(A)\to \mathcal{A}$.
    \item
    The map $H(A)\to H(\Gamma)$ is injective, and identifies half-edges $H(v)$ of $A$ with legs of $\Gamma_v$.
    \item
    The maps $V(\Gamma)\to V(A)$ and $H(A)\to H(\Gamma)$ respect the incidence relation of half-edges and vertices, and the pairs of half-edges forming edges. (In particular, $E(A)\subset E(\Gamma)$.)
\end{enumerate}
\end{definition}
Given an $A$-structure $\Gamma\to A$ on $\Gamma$, the map $m_\Gamma: H(\Gamma)\to \mathbb{Z}_{\geq 0}$ and the injection $H(A)\to H(\Gamma)$ defines the map $$m_A: H(A)\to \mathbb{Z}_{\geq 0}.$$
Then there is a natural gluing map
\begin{equation*}
    \xi_{(\Gamma\to A, m_\Gamma, a)}: \mathfrak{M}_{(\Gamma, m_\Gamma, a)}^{tw, triv}\to \mathfrak{M}_{(A, m_A, a)}^{tw, triv}.
\end{equation*}
We have
\begin{enumerate}
    \item 
    If $H(A)\to H(\Gamma)$ maps $i\in H(v)\subset H(A)$ to $j\in H(w)\subset H(\Gamma)$, then $\xi_{(\Gamma\to A, m_\Gamma, a)}^*\psi_{v,i}=\psi_{w, j}$.
    \item
    $\xi_{(\Gamma\to A, m_\Gamma, a)}^*\kappa_{v,l}=\sum_{w\in V(\Gamma_v)} \kappa_{w, l}$.
\end{enumerate}
These can be proven by ways analogous to their counterparts for moduli stacks of prestable curves, see \cite[Section 3.2]{bs}. These properties allow, in a straightforward manner, a description of the pull-back $\xi_{(\Gamma\to A, m_\Gamma, a)}^*\alpha$ of a decoration $\alpha$ on $\mathfrak{M}_{(A,m_A, a)}^{tw, triv}$ as in (\ref{eqn:deco}).

Suppose a prestable graph $\Gamma$ admits $A$ and $B$-structures $f_A: \Gamma\to A$ and $f_B:\Gamma\to B$. As in \cite[Section 3.2]{bs}, the pair $(f_A, f_B)$ is called a generic $(A,B)$-structure on $\Gamma$ if each half-edge of $\Gamma$ corresponds to either a half-edge of $A$ or a half-edge of $B$. In other words, the injections $H(A)\to H(\Gamma)$ and $H(B)\to H(\Gamma)$ have disjoint images and their union is the whole $H(\Gamma)$. 

With a generic $(A, B)$ structure on $\Gamma$, the maps $m_A: H(A)\to \mathbb{Z}_{\geq 0}$ and $m_B: H(B)\to \mathbb{Z}_{\geq 0}$ uniquely determine a map $m_\Gamma: H(\Gamma)\to \mathbb{Z}_{\geq 0}$.

\begin{proposition}
Let $A$ and $B$ be $\mathcal{A}$-valued prestable graphs for $\mathfrak{M}_{g, I, m, a}$. Let $m_A$ and $m_B$ be as in (\ref{eqn:orb_str}). Then the fiber product of the gluing maps $\xi_{(A, m_A, a)}: \mathfrak{M}_{(A,m_A,a)}^{tw, triv}\to \mathfrak{M}_{g, I,m,a}^{tw, triv}$ and $\xi_{(B, m_B, a)}: \mathfrak{M}_{(B,m_B,a)}^{tw, triv}\to \mathfrak{M}_{g, I,m,a}^{tw, triv}$ is given by the disjoint union
\begin{equation*}
    \coprod_{\Gamma\in \mathcal{G}_{A,B}}\mathfrak{M}_{(\Gamma, m_\Gamma, a)}^{tw, triv}
\end{equation*}
taken over the set $\mathcal{G}_{A,B}$ of isomorphism classes of generic $(A,B)$ structures on prestable graphs $\Gamma$, with two projection maps given by $\xi_{(\Gamma\to A, m_\Gamma, a)}$ and $\xi_{(\Gamma\to B, m_\Gamma, a)}$. Moreover, the top Chern class of the excess bundle 
\begin{equation*}
    E_{(\Gamma, m_\Gamma, a)}=\xi_{(\Gamma\to A, m_\Gamma, a)}^*\mathcal{N}_{\xi_{(A, m_A, a)}}/\mathcal{N}_{\xi_{(\Gamma\to B, m_\Gamma, a)}}
\end{equation*}
is given by 
\begin{equation}\label{eqn:ctop}
    c_{top}(E_{(\Gamma, m_\Gamma, a)})=\prod_{e=(h, h')\in E(A)\cap E(B)\subset E(\Gamma)} (-\psi_h-\psi_{h'})/m_\Gamma(h).
\end{equation}
\end{proposition}
The proof of this Proposition is similar to that of \cite[Proposition 3.6]{bs}. The only thing to notice is the factor $1/m_\Gamma(h)$ in (\ref{eqn:ctop}), which accounts for the fact that the gerbe at nodes are not trivialized. 

The projection formula yields the following description of products of decorated stratum classes:
\begin{equation}\label{eqn:product}
[(A, m_A, a), \alpha]\cdot [(B, m_B, a), \beta]=\sum_{\Gamma\in \mathcal{G}_{A, B}}(\xi_{(\Gamma, m_\Gamma, a)})_*(\xi_{(\Gamma\to A, m_\Gamma, a)}^*\alpha \cdot \xi_{(\Gamma\to B, m_\Gamma, a)}^*\beta\cdot c_{top}(E_{(\Gamma, m_\Gamma, a)}).
\end{equation}
The above discussion shows that the right-hand side of (\ref{eqn:product}) is contained in the span of decorated stratum classes. This shows that the span is closed under intersection product.

\subsection{Closure under pushforwards}
\subsubsection{Gluing maps}
Pushing forward via (\ref{eqn:gluing_map}) is straightforward to understand. Let $\Gamma_0$ be a graph, $m_{\Gamma_0}: H(\Gamma_0)\to \mathbb{Z}_{\geq 0}$, and $a_0: V(\Gamma_0)\to \mathcal{A}$ an $\mathcal{A}$-valuation. Suppose for each $v\in V(\Gamma_0)$, we have the data of a genus $g(v)$ graph $\Gamma_v$, a map $m_{\Gamma_v}:H(\Gamma_v)\to \mathbb{Z}_{\geq 0}$, and an $\mathcal{A}$-valuation $a_v: V(\Gamma_v)\to \mathcal{A}$, satisfying the following conditions:
\begin{enumerate}
    \item For each $v\in V(\Gamma_0)$, there is an embedding of the set of half-edges $H(\Gamma_0)|_v\subset H(\Gamma_v)$ so that $m_{\Gamma_v}$ restricts to $m_{\Gamma_0}|_{H(\Gamma_0)|_v}$.
    \item For each $v\in V(\Gamma_0)$, $\sum_{v'\in V(\Gamma_v)}a_v(v')=a_0(v)$.
\end{enumerate}
For each $v\in V(\Gamma_0)$, we can replace $v$ by the graph $\Gamma_v$ and replace half-edges in $H(\Gamma_0)|_v$ by their images in $H(\Gamma_v)$. This yields a new graph $\Gamma$ with $V(\Gamma)=\cup_{v\in V(\Gamma_0)}V(\Gamma_v)$ and $H(\Gamma)=\cup_{v\in V(\Gamma_0)}H(\Gamma_v)$. Define $m_\Gamma$ by $m_\Gamma|_{H(\Gamma_v)}=m_{\Gamma_v}$ and define $a:V(\Gamma)\to \mathcal{A}$ by $a|_{V(\Gamma_v)}=a_v$.

This construction describes a natural gluing map
\begin{equation*}
    \prod_{v\in V(\Gamma_0)}\mathfrak{M}_{(\Gamma_v, m_{\Gamma_v}, a_v)}^{tw, triv}\to \mathfrak{M}_{(\Gamma, m_\Gamma, a)}^{tw, triv}.
\end{equation*}
The pushforward of a class $\prod_{v\in V(\Gamma_0)}[(\Gamma_v,m_{\Gamma_v}, a_v), \alpha_v]$ under this gluing map is given by $$\frac{1}{\prod_{h\in E(\Gamma_0)}m_{\Gamma_0}(h)}[(\Gamma, m_\Gamma, a), \alpha]$$ where $\alpha$ is formed by combining $\alpha_v$ using $V(\Gamma)=\cup_{v\in V(\Gamma_0)}V(\Gamma_v)$.

\subsubsection{Forgetful maps}
We begin with an analogue of \cite[Proposition 3.11]{bs}. 

Let $[(\Gamma, m_\Gamma, a), \alpha]\in R^*(\Mfrak_{g,I\cup\{\bullet\}, m',a}^{tw, triv})$ with $\alpha=\prod_{v\in V(\Gamma)}\alpha_v$. Let $v\in V(\Gamma)$ be the vertex incident to $\bullet$ and let $\Gamma'$ be the graph obtained from $\Gamma$ by forgetting $\bullet$ and stabilizing if the vertex $v$ becomes unstable. $v$ become unstable exactly when $g(v)=0$, $a(v)=\mathbb{0}$, and $v$ has two other half-edges. By stability, the curve component corresponding to $v$ is a genus $0$ stable twisted curve with three marked points. By assumption the marked point corresponding to $\bullet$ has trivial stack structure: $m_\Gamma(\bullet)=1$. Stabilization contracts this curve component, hence the other marked points have opposite stack structure. Hence $m_{\Gamma'}:H(\Gamma')\to \mathbb{Z}_{\geq 0}$ is well-defined from $m_\Gamma$. $a':V(\Gamma')\to \mathcal{A}$ is naturally defined from $a$.
\begin{proposition}\label{prop:forget_bullet}
\hfill
\begin{enumerate}
    \item If $v$ remains stable after forgetting $\bullet$, then 
    \begin{equation*}
        (\pi_\bullet)_*[(\Gamma, m_\Gamma, a), \alpha]=(\xi_{(\Gamma',, m_{\Gamma'}, a')})_*\left((\pi_v)_*\alpha_v\cdot \prod_{w\neq v}\alpha_w\right),
    \end{equation*}
     where $\pi_v$ is the forgetful map of $\bullet$ of the vertex $v$.
    \item Suppose $v$ becomes unstable after forgetting $\bullet$. If $\alpha_v\neq 1$, then $[(\Gamma, m_\Gamma, a), \alpha]=0$. If $\alpha_v=1$, then \begin{equation*}
      (\pi_\bullet)_*[(\Gamma, m_\Gamma, a), \alpha]=[(\Gamma', m_{\Gamma'}, a'),\prod_{w\neq v}\alpha_w].  
    \end{equation*}
\end{enumerate}
\end{proposition}
The proof of Proposition \ref{prop:forget_bullet} is similar to that of \cite[Proposition 3.11]{bs} and is omitted.

We can study pushforwards of products of $\kappa$ and $\psi$ classes under $\pi_\bullet$, as follows. A product 
\begin{equation*}
    \prod_b \kappa_b^{e_b}\cdot\prod_{i\in I} \psi_i^{l_i}\cdot \psi_{\bullet}^{l_\bullet}\in R^*(\Mfrak_{g,I\cup\{\bullet\},m', a}^{tw, triv})
\end{equation*}
can be written as 
\begin{equation*}
    (\pi_\bullet)^*(\prod_b \kappa_b^{e_b}\cdot\prod_{i\in I} \psi_i^{l_i})\cdot \psi_{\bullet}^{l_\bullet} + \text{boundary terms}.
\end{equation*}
This uses a straightforward analogue of \cite[Proposition 3.10]{bs} and known intersection formulas on $\Mfrak_{g,I\cup\{\bullet\},m', a}^{tw, triv}$. Applying $(\pi_\bullet)_*$ yields
\begin{equation*}
    \left(\prod_b \kappa_b^{e_b}\cdot\prod_{i\in I} \psi_i^{l_i}\right)\cdot \kappa_{l_\bullet-1}+(\pi_\bullet)_*(\text{boundary terms}).
\end{equation*}
The boundary terms can be handled by induction together with Proposition \ref{prop:forget_bullet}.

We conclude that $(\pi_\bullet)_*$ preserves the span of decorated stratum classes.

\end{document}